# Notes on the Golden Ratio: The Golden Rule of Vector Similarities in Space

Artyom M. Grigoryan*[1] and Meruzhan M. Grigoryan

[1] Department of Electrical and Computer Engineering, The University of Texas at San Antonio, USA
and The Yerevan State University, Armenia
amgrigoryan@utsa.edu

**Abstract -** In this work, we have abstractly generalized the similarity law for multidimensional vectors. Initially, the law of similarity was derived for one-dimensional vectors. Although it operated with such values of the ratio of parts of the whole, it meant linear dimensions (a line is one-dimensionality). The concept of the general golden ratio (GGR) for the vectors in the multidimensional space is presented and described in detail with equations and solutions. It shown that the GGR depends on the angles. Main properties of the GGR are given with illustrative examples. We introduce and discuss the concept of the similar vectors and the set of similarities for a given vector. Also, we present our vision on the theory of the golden ratio for triangles and describe the similar triangles with examples.

## 1. Introduction

The golden ratio, or the famous number $\Phi = 1.6180\ldots$, is known as the divine proportion between two quantities [1,2]. For positive numbers $a$ and $b$, this ratio equals to $\Phi = b/a = (b+a)/b$, if $b > a$. The number $\Phi$ was found in proportions of parts in construction of Egyptian pyramids [3], in the art of painting [4]-[6], in fashion [7], in medicine [8,9], face detection [10], image enhancement [11,12], and many applications in engineering.

In this work we generalize and describe the Golden ratio in the multi-dimensional vector space. The concept of the general golden ration (GGR) is described in detail with examples in 2, 3, and 6-dimensional spaces. We present our vision on this concept, the golden ratio that depends on the direction in the space. Our vision on the problem of this divine proportion can be formulated as follows:

1. The concept of the golden ratio has generalization in the multi-dimensional space;
2. The general golden ratio is function of one or a few angles. It is one of four solutions of the golden equation presented here;
3. To each vector $V$ in the space corresponds a set of similarity, or a set of vectors being in golden ratio with $V$;
4. GGR can be used to describe similar figures. As an example, the vector space of triangles is described;
5. Each vector affects its environment, stimulating its influence through the imposition of its likeness.

The rest of the paper is organized in the following way. Section 2 presents the definition of general golden ratio in the space with examples. The main equation of the GGR is described in Section 3. The golden ratio in the 2-D space with examples is considered in Section 4. In Section 5, the vector similarity sets are described and the 3-D space is considered. The vector space of triangles and golden ratio of triangles are described in Section 6. Illustrative examples are given.

## 2. The Generalized Golden Ratio (GGR)

In this section, we extend the well-known concept of the gold ratio in the $n$-dimension vector space. Let $V$ be the normed vector space over real numbers $\lambda$. It means that exists such a real-valued function $l: V \to R$, which is denoted by $l(x) = \|x\|$ and has the following properties:

1. $\|x\| \geq 0$, and $\|x\| = 0$ means that $x = 0$;
2. $\|\lambda x\| = \lambda \|x\|$, for any real number $\lambda \in R$;
3. $\|x + y\| \leq \|x\| + \|y\|$, for any elements $x$ and $y \in V$.

The function $\|x\|$ is called the length, or the norm, of the element $x$.

***Definition 1.*** In a normed vector space $V$, the function $f: V \times V \to R$ which satisfies the conditions

1. $f(a, b) \geq 0$, for any elements $a$ and $b \in V$;
2. If $\|a\| = 0$, $f(a, b) = 0$;
3. If $\|b\| \to 0$, $f(a, b) \to \infty$;

4. For any real number $\lambda \in R$, $f(\lambda a, \lambda b) = f(a, b)$,

is called *the proportion*.
As an example, when $V = R$, the function $f(a, b) = |a|/|b|$ is the proportion.

***Definition 2.*** Given a proportion $f: V \times V \to R$ in the normed vector space $V$, two elements $a$ and $b \in V$ are called *the golden pair*, if the following holds:

$$f(b + a, b) = f(a, b). \tag{1}$$

The elements $a$ and $b$ are also called *the golden ratio elements*, or *the golden pair*.

***Definition 3.*** If the elements $a$ and $b$ are golden ratio elements, for a proportion $f$, then the ratio $\Phi = \|a\|/\|b\|$ is called *the generalized golden ratio*, or shortly GGR.

In the example below for the 1-D vector space, this number is the known golden ratio [1]. Therefore, we use a similar name. In order not in any way to underestimate the historical nature of this number, we added a generalized meaning.

***Example 1:*** Let $V = R$ and the function $f(a, b) = |a|/|b|$ be the proportion. Here, we consider that $\|a\| = |a|$ for real numbers. To find the golden ratio elements $a$ and $b$, we consider the primary equation

$$\frac{|b + a|}{|b|} = \frac{|b|}{|a|}, \quad \text{or} \quad \frac{|b + a|}{|a|} = \frac{|b|^2}{|a|^2}.$$

It can be written as

$$\left(\frac{|b|}{|a|}\right)^2 = \frac{|b + a|}{|a|} = \left|\frac{|b|}{|a|} \frac{b}{|b|} + \frac{a}{|a|}\right|,$$

or

$$\Phi^2 = |\Phi + \text{sign}(ab)|. \tag{2}$$

The solutions of this equation are considered for the following two possible cases.

1. Case $\Phi + \text{sign}(ab) \geq 0$: Then, the equation to be solved is

$$\Phi^2 - \Phi - \text{sign}(ab) = 0. \tag{3}$$

The solutions are

$$\Phi_{1,2} = \frac{1 \pm \sqrt{1 + 4\,\text{sign}(ab)}}{2}$$

Such numbers are real only when $\text{sign}(ab) = 1$. Therefore,

$$\Phi_1 = \frac{1 + \sqrt{5}}{2} = 1.6180\ldots, \quad \Phi_2 = \frac{1 - \sqrt{5}}{2} = -0.6180\ldots. \tag{4}$$

A golden ratio is a positive number. The first number $\Phi_1$ is considered, but the second $\Phi_2$ is not, since it is negative. Thus, the golden ratio is $\Phi = \Phi_1$.

2. Case $\Phi + \text{sign}(ab) < 0$: The following equation is considered:

$$\Phi^2 + \Phi + \text{sign}(ab) = 0 \tag{5}$$

with the solutions

$$\Phi_{1,2} = \frac{1 \pm \sqrt{1 - 4\,\text{sign}(ab)}}{2}.$$

Such numbers are real only when $\text{sign}(ab) = -1$. Therefore,

$$\Phi_1 = \frac{1+\sqrt{5}}{2} = 1.6180\ldots, \qquad \Phi_2 = \frac{1-\sqrt{5}}{2} = -0.6180\ldots.$$

Here, $\Phi_2$ is negative and $\Phi_1$ also needs to be discarded, since the condition of consideration is violated, that is, $\Phi_1 + \text{sign}(ab) = \Phi_1 - 1 = 0.6180 > 0$. Thus, a golden pair cannot be composed from numbers of opposite signs.

This example shows that for any number $a \neq 0$ there is only one number that composes the golden ratio with $a$. This number is $a\Phi_1$. The golden pair is $\{a, a\Phi_1\}$.

***Example 2:*** Consider the $n$-dimensional vector space $V = R^n, n > 1$. It is not difficult to show that the function

$$f(a,b) = \frac{\|a\|}{\|b\|}, \quad a,b \in V, \quad b \neq 0, \tag{6}$$

is the proportion. Here, the norms $\|a\| = \sqrt{a_1^2 + a_2^2 + \cdots + a_n^2}$ and $\|b\| = \sqrt{b_1^2 + b_2^2 + \cdots + b_n^2}$.

We consider the golden ratio equation (rule) $f(b+a, b) = f(a, b)$ written as

$$\frac{\|b+a\|}{\|b\|} = \frac{\|b\|}{\|a\|}, \quad \text{or} \quad \|b+a\|\|a\| = \|b\|^2.$$

The following calculations are valid for this equation:

$$\sqrt{\|a\|^2 + \|b\|^2 + 2(a,b)} \, \|a\| = \|b\|^2,$$
$$\sqrt{\|a\|^2 + \|b\|^2 + 2\|a\|\|b\|\cos(\alpha)} \, \|a\| = \|b\|^2,$$
$$1 + \frac{\|b\|^2}{\|a\|^2} + 2\frac{\|b\|}{\|a\|}\cos(\alpha) = \frac{\|b\|^4}{\|a\|^4}.$$

Here, $(a,b)$ is the inner product in the space $V$, and $\alpha$ is the angle between the vectors $a$ and $b$. Denoting the golden ratio $x = \Phi = \|b\|/\|a\|$, we obtain the equation $1 + x^2 + 2x\cos(\alpha) = x^4$, or

$$x^4 - x^2 - 2x\cos(\alpha) - 1 = 0. \tag{7}$$

This equation has four roots $x_n = x_n(\alpha)$, $n = 1:4$, which are functions of the angle. Thus, the GGR depends on the angle. Since $\cos(\alpha) = \cos(2\pi - \alpha)$, the roots $x_n(2\pi - \alpha) = x_n(\alpha)$, for $\alpha \in [0, \pi]$. Also, $x_n(\pi + \alpha) = -x_n(\alpha)$.

1. Case $\alpha = 0$ (vectors are collinear, or are in the same directions): The equation $x^4 - x^2 - 2x - 1 = 0$ can be written as

$$x^4 - x^2 - 2x - 1 = (x^2 - x - 1)(x^2 + x + 1) = 0$$

Therefore, we consider only the solutions of the equation $x^2 - x - 1 = 0$, which are

$$x_{1,2} = \frac{1 \pm \sqrt{5}}{2} = \Phi_{1,2}.$$

The positive solution is $x_1 = \Phi_1 = 1.6180$.

2. Case $\alpha = \pi/2$ (vectors are perpendicular): The equation $x^4 - x^2 - 1 = 0$ is reduced to two equations

$$x^2 = \frac{1+\sqrt{5}}{2}, \frac{1-\sqrt{5}}{2}.$$

Therefore, the first equation is considered and its two solutions are

$$x_{1,2} = \pm\sqrt{\frac{1+\sqrt{5}}{2}} = \pm\sqrt{\Phi_1} = \pm\sqrt{1.6180}. \tag{8}$$

The positive solution is $x_1 = \sqrt{\Phi_1} = 1.2720$.

3. Case $\alpha = \pi$ (vectors are collinear in the opposite directions): The equation $x^4 - x^2 + 2x - 1 = 0$ can be written as

$$x^4 - x^2 + 2x - 1 = (x^2 + x - 1)(x^2 - x + 1) = 0$$

Therefore, we consider the solution of equation $x^2 + x - 1 = 0$, which are

$$x_1 = \frac{-1 - \sqrt{5}}{2} = -\Phi_1, \quad x_2 = \frac{-1 + \sqrt{5}}{2} = -\Phi_2.$$

The positive solution is $x_2 = -\Phi_2 = 0.6180$.

4. When the ratio $x = 1$, Eq. 7 is $2\cos(\alpha) + 1 = 0$. Thus, $\alpha = 180 \pm 60$. Two vectors with the angle $\alpha = 240$ or 120 (in degree) between them compose a golden pair when their lengths are equal. Figure 1 shows three vectors $a_1, a_2$, and $a_3$ with the same length. Each of these vectors is the golden pair of other two.

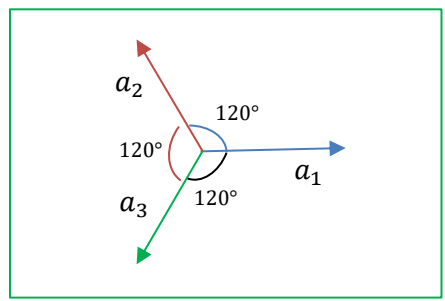

**Fig. 1** The golden pairs $\{a_1, a_2\}$, $\{a_1, a_3\}$, and $\{a_2, a_3\}$.

## 3. Main Equation of Golden Ration and Its Analytical Solution

Consider the positive root of Eq. 7, $x^4 - x^2 - 2x\cos(\alpha) - 1 = 0$, where $x$ is a function of the angle $x = x(\alpha)$. The equation is

$$x^4(a) - x^2(a) - 2\cos\alpha\, x(\alpha) - 1 = 0 \tag{9}$$

with initial condition $x(0)=1.6180...$ . For each angle $\alpha$, the polynomial of order 4 in Eq. 9 has four roots, and two of them are complex and therefore complex conjugate. We are looking for a positive solution of this equation, which should be only one. Equation 9 can be written as $x^4 - x^2 - 1 = 2\cos(\alpha)x$. The parabola $P(x) = (x^2)^2 - x^2 - 1$ crosses the straight line with the slope $2\cos(\alpha)$ at two points. As illustration, Fig. 2 shows the graph of the polynomial $P(x)$ together with the line $2\cos(\alpha)x$, when $\alpha = 45°$ in part (a) and $\alpha = 100°$ in part (b).

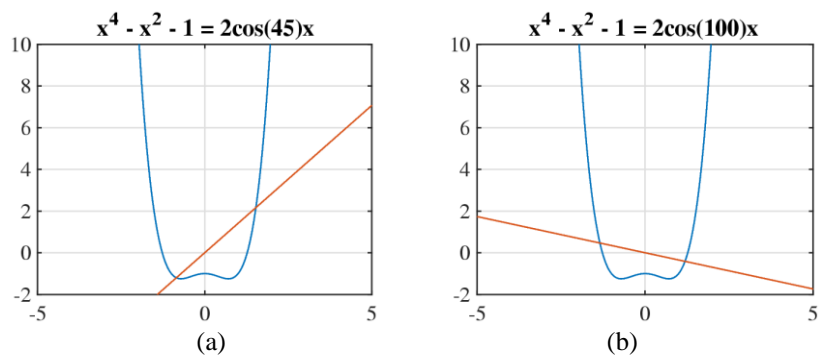

**Fig. 2** The graphs of the polynomial and line, when (a) $\alpha = 45°$ and (b) $\alpha = 100°$.

For the angle $\alpha = 45°$, the solutions of Eq. 9 are

$$x_1 = -0.8492, \quad x_2 = 1.5322, \quad x_3 = -0.3416 - 0.8072i, \quad x_4 = \bar{x}_3 = -0.3416 + 0.8072i.$$

The numbers are written here with 4 decimal precision. At points $x_1$ and $x_2$, the parabola $P(x) = (x^2)^2 - x^2 - 1$ crosses the straight line $y = \sqrt{2}x$, as shown in part (a). The second coordinate is positive. Thus, the required root $\Phi = x(45°) = 1.5322$.

For the angle $\alpha = 100°$, the solutions of Eq. 9 are

$$x_1 = -1.3455, \quad x_2 = 1.1894, \quad x_3 = 0.0780 + 0.7865i, \quad x_4 = \bar{x}_3 = 0.0780 - 0.7865i.$$

The parabola $P(x)$ crosses the straight line $y = 2\cos(100)x$ at the point $x_1 < 0$ and the positive point $x_2$. Therefore, $\Phi = x(100°) = 1.1894$.

## 1. Similarity equation and its roots

Consider again the main equation

$$x^4 - x^2 - 2x\cos(\alpha) - 1 = 0. \tag{7}$$

We call the continuous roots of this equation the similarity functions and denote them by the symbols $x_k(\alpha), k \in 1:4$. To analyze this equation of the generalized golden ratio, we consider its roots. Figure 3 shows the graphs of four roots $x_n(\alpha), n = 1:4$, of the equation. The real and imaginary parts of the roots are shown blue and red colors, respectively. The angles $\alpha$ are in the interval $[0, 2\pi]$ with step 0.0015 radians, or 1/12 degrees. These roots were calculated by using MATLAB function 'roots.m' with the command 'x=roots([1,0,-1,-2*cos(a),-1])'.

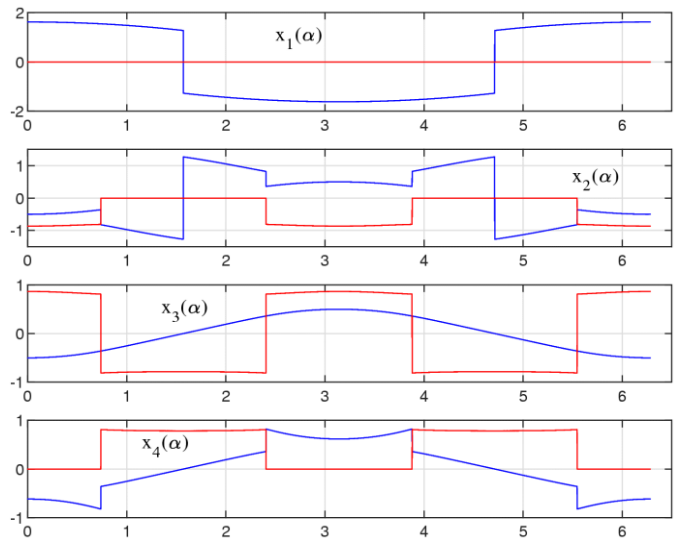

**Fig. 3** Graphs of the four roots of Eq. 7.

The graphs of the roots are symmetric with respect to the vertical at angle-point $\alpha = \pi$. In some parts these functions change sign. For example, the change of sign in the real part of the first solution $x_1(\alpha)$ occurs at angles $\pi/2$ and $3\pi/2$ and the jump is equals to $2\sqrt{\Phi_1} = 2 \times 1.2720$. For other roots, the discontinuities can be seen at angle-points $\pi/2 \pm \pi/4$ and $3\pi/2 \pm \pi/2$.

In Figure 4, these four roots are plot in polar form. The first plot is like the apple, the 2$^{nd}$ as a four-petal flower, the 3$^{rd}$ as the egg (Earth), and the 4$^{th}$ plot is an unknown figure for us.

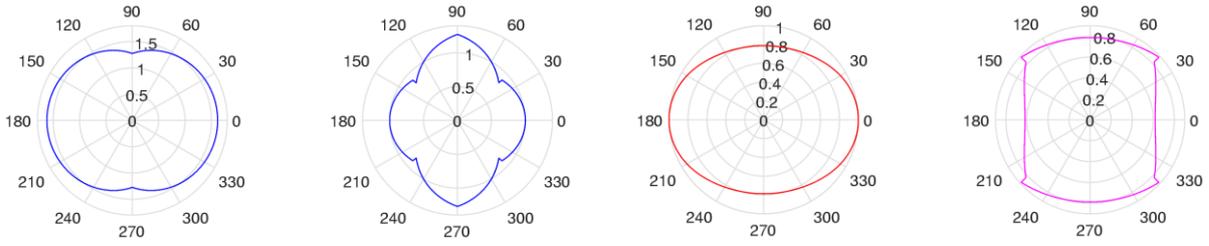

**Fig. 4** Polar plot of the four roots of Eq. 7.

It is not difficult to note that the following holds for the roots of the above equation: $x_1(\alpha) + x_2(\alpha) + x_3(\alpha) + x_4(\alpha) = 0$. Thus, $x_4(\alpha)$ equals to the sum of the first three roots with sign minus. The 4$^{\text{th}}$ plot is the sum $x_1(\alpha) + x_2(\alpha) + x_3(\alpha)$ in polar form. Figure 5 shows the plots of the sum of roots $x_1(\alpha) + x_2(\alpha)$, $x_2(\alpha) + x_3(\alpha)$, and $x_1(\alpha) + x_3(\alpha)$ in polar form in parts (a), (b), and (c), respectively.

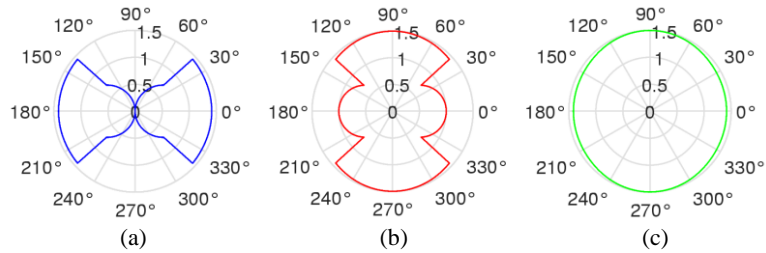

(a)      (b)      (c)

**Fig. 5** The polar plots of the sum of two roots of Eq. 7.

## 2. Analyze of Solutions

We will regroup the obtained set of roots $x_1(\alpha)$, $x_2(\alpha)$, $x_3(\alpha)$, and $x_4(\alpha)$ of the above equation in the following way. It is not difficult to note from Fig. 3 that, for each angle, there are two real solutions of Eq. 7. Even more, there is only one positive solution for each angle. Figure 6(a) shows these two roots (solutions), which we denote by $\Phi_1(\alpha)$ and $\Phi_2(\alpha)$. These solutions are composed as follows:

$$\Phi_1(\alpha) = \begin{cases} x_1(\alpha), & \text{if } x_1(\alpha) \geq 0; \\ x_2(\alpha), & \text{if } x_2(\alpha) \geq 0, \end{cases} \quad \Phi_2(\alpha) = \begin{cases} x_2(\alpha), & \text{if } x_1(\alpha) \geq 0; \\ x_1(\alpha), & \text{if } x_2(\alpha) \geq 0, \end{cases} \tag{10}$$

In part (b), the graph of the sum of these solutions is shown, $\Psi(\alpha) = \Phi_1(\alpha) + \Phi_2(\alpha)$. One can note that $|\Psi(\alpha)| \leq 1$.

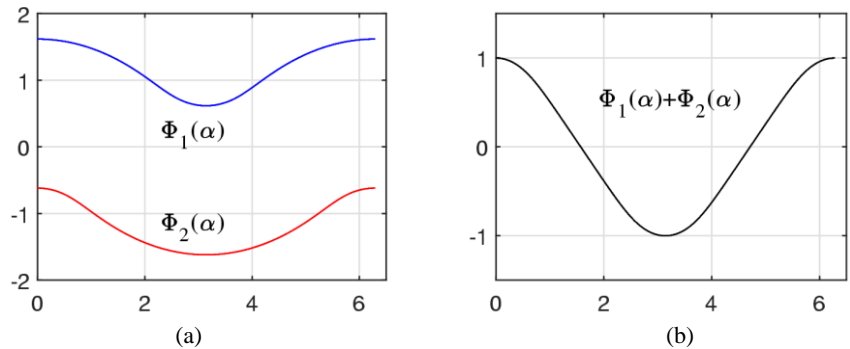

(a)      (b)

**Fig. 6** Two real solutions of Eq. 7 and their sum.

It is interesting to note that $\Phi_1(\alpha + \pi) = -\Phi_2(\alpha), \alpha \in [0, 2\pi]$.

The MATLAB-based codes for computing the solutions $\Phi_k(\alpha), k = 1,2$, are given below.

```
% call: golden_ratio1.m
% 1st real root of the equation

function r = golden_ratio1(a)
    f0=2*cos(a);
    r=roots([1,0,-1,-f0,-1]);
    kk=find(imag(r)==0);
    r2=r(kk);
    if r2(1)>0
        r=r2(1);
    else
        r=r2(2);
    end
end
```

```
% call: golden_ratio2.m
% 2nd real root of the equation

function r=golden_ratio2(a)
    f0=2*cos(a);
    r=roots([1,0,-1,-f0,-1]);
    kk=find(imag(r)==0);
    r2=r(kk);
    if r2(1)>0
        r=r2(2);
    else
        r=r2(1);
    end
end
```

Figure 7 shows the graph of the positive roots $\Phi_1(\alpha)$ calculated in the interval of interval $[0,2\pi]$. We call the function $\Phi(\alpha) = \Phi_1(\alpha)$ with this graph *the general golden ratio function*, or *the GGR function*. For this function, the minimum is 0.6180 at the angle-point $\pi$, and the maximum is 1.6180 at 0 and $2\pi$. The Golden ratio equals to 1 at angles $2/3\pi$ and $4/3\pi$ (as shown in Fig. 7). The mean of the Golden ratios in this interval equals to 1.192880, approximately at angles 1.7385 and 4.5447, or 99.6075 ° and 260.3925° (or 180° ∓ 80.3925°).

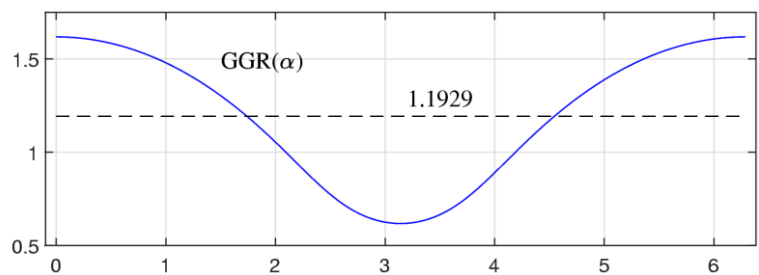

**Fig. 7** The Golden ratio function $\Phi_1(\alpha), \alpha \in [0,2\pi]$.

The GGR function $\Phi(\alpha)$ has approximately the form of the cosine function. Together with the GGR function, the following function is shown in Fig. 8:

$$y(\alpha) = \frac{1}{2}[\Phi_1(0) + \Phi_1(\pi) + \cos\alpha] = \frac{1}{2}[1.6180 + 0.6180] + \frac{1}{2}\cos\alpha = 1.1180 + \frac{1}{2}\cos\alpha. \quad (11)$$

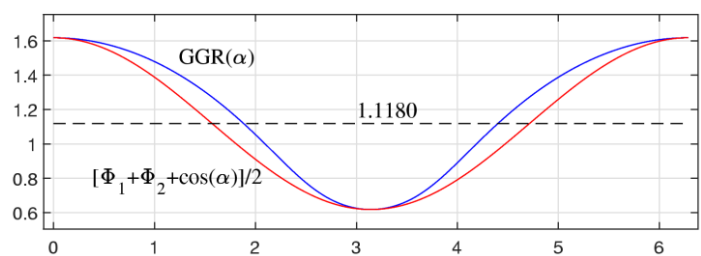

**Fig. 8** The General Golden Ratio function.

Figure 9 shows the graph of the GGR function versus angles in degrees. A few points on the graph are marked for the values of this function at angles 36°, 72°, 108°, and 144°, plus at the angle 290.70°, at which the Golden ratio equals to $\sqrt{2}$.

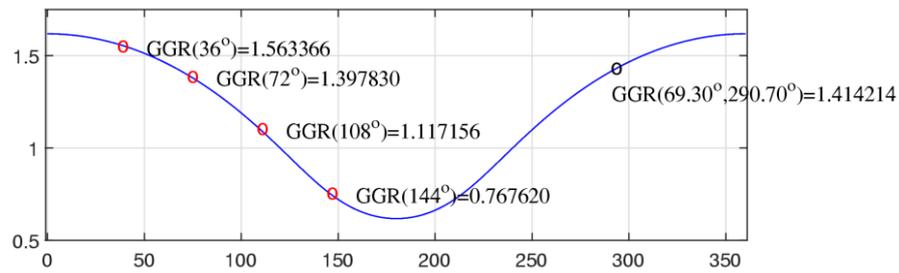

**Fig. 9** The General Golden ratio function $\Phi_1(\alpha), \alpha \in [0, 2\pi]$, with a few values on it.

The complex roots of Eq. 7, can also be regrouped, by using the phases of two complex solutions $x_3(\alpha) = |x_3(\alpha)|e^{i\vartheta_3(\alpha)}$ and $x_4(\alpha) = |x_4(\alpha)|e^{i\vartheta_4(\alpha)}$. Namely, these functions are calculated as

$$\Phi_3(\alpha) = \begin{cases} x_3(\alpha), & \text{if } \vartheta_3(\alpha) > 0; \\ x_4(\alpha), & \text{if } \vartheta_4(\alpha) > 0, \end{cases} \quad \Phi_4(\alpha) = \begin{cases} x_4(\alpha), & \text{if } \vartheta_3(\alpha) > 0; \\ x_3(\alpha), & \text{if } \vartheta_4(\alpha) > 0. \end{cases} \quad (12)$$

Note that $\vartheta_4(\alpha) = -\vartheta_3(\alpha)$.

MATLAB-based codes for the complex roots of Eq. 7:

```
% call: golden_ratio3.m
% 1st complex root of the equation

function r=golden_ratio3(a)

    f0=2*cos(a);
    r=roots([1,0,-1,-f0,-1]);
    kk=find(abs(imag(r))>0);
    r3=r(kk);
    if phase(r3(1))>0
        r=r3(1);
    else
        r=r3(2);
    end
end
```

```
% call: golden_ratio4.m
% 2nd complex root of the equation

function r=golden_ratio4(a)

    f0=2*cos(a);
    r=roots([1,0,-1,-f0,-1]);
    kk=find(abs(imag(r))>0);
    r3=r(kk);
    if phase(r3(1))>0
        r=r3(2);
    else
        r=r3(1);
    end
end
```

Figure 10 shows the graphs of the real and imaginary parts of the complex solution $\Phi_3(\alpha)$. One can note that the real part of the solution $\Phi_3(\alpha)$ has values in the interval $[-0.5, 0.5]$.

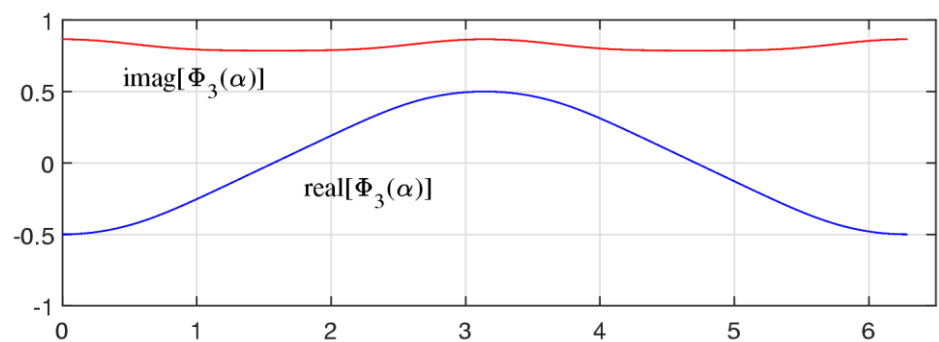

**Fig. 10** The complex solution $\Phi_3(\alpha), \alpha \in [0, 2\pi]$, of Eq. 7.

**Properties of the roots**

The following properties hold for the roots of the equation:

$$\Phi_k(-\alpha) = \Phi_k(\alpha) \text{ and } \Phi_k(\alpha + 2\pi) = \Phi_k(\alpha), k = 1,2,3,4, \tag{13}$$

$$\Phi_1(\alpha) > 0 \text{ and } \Phi_2(\alpha) < 0, \text{ for all } \alpha, \tag{14}$$

$$\Phi_3(\alpha), \Phi_4(\alpha) \in C \text{ and } \Phi_3(\alpha) = \overline{\Phi_4}(\alpha), \tag{15}$$

Also, from Eq. 7, we obtain the following identities:

$$\Phi_1(\alpha) + \Phi_2(\alpha) + \Phi_3(\alpha) + \Phi_4(\alpha) = 0 \tag{16}$$

$$\Phi_1(\alpha)\Phi_2(\alpha)\Phi_3(\alpha)\Phi_4(\alpha) = -1, \tag{17}$$

$$\sum_{i \neq j}^{4} \Phi_i(\alpha)\Phi_j(\alpha) = -1, \tag{18}$$

$$\sum_{i \neq j \neq k}^{4} \Phi_i(\alpha)\Phi_j(\alpha)\Phi_k(\alpha) = -2\cos(\alpha). \tag{19}$$

Due to Eqs. 15 and 16, the real part $\Re(\Phi_3(\alpha))$ of the 3rd solution shown in Fig. 10 equals to $-[\Phi_1(\alpha) + \Phi_2(\alpha)]/2$ (see Fig. 6(b)).

It is also not difficult to see that the solutions are transformed into each other under the transformation $\alpha \to \alpha + \pi$. Indeed, the following identities are valid, for any angle $\alpha$:

$$\Phi_1(\alpha) = -\Phi_2(\alpha + \pi), \quad \Phi_3(\alpha) = -\Phi_4(\alpha + \pi). \tag{20}$$

## 4. Examples of Golden Ratios

*Example 1* (1-D vectors):

For 1-D vectors (real numbers), or the elements of the real line $R$, we define the inner product as $(a, b) = ab$. Then, the angle is defined as

$$\cos(\alpha) = \frac{ab}{\|a\|\|b\|} = \frac{ab}{|a||b|} = \text{sign}(ab),$$

which means that the angle between similar elements may take only two values, 0 and $\pi$. Therefore, the set of similarity, that is, the set of numbers that are in golden ratios with the number $a$, is defined as

$$S(a) = \{a\} = \{|a|\Phi(\arccos(\text{sign}(ab)))\text{sign}(b); b = \pm 1\}. \tag{21}$$

When $a > 0$, that is, $\text{sign}(a) = 1$, the number $a$ is in golden ratio with numbers of the set

$$S(a) = \{|a|\Phi(0), -|a|\Phi(\pi)\}.$$

When $a < 0$, that is, $\text{sign}(a) = -1$, the golden pairs are defined by the similar set

$$S(a) = \{|a|\Phi(\pi), -|a|\Phi(0)\}.$$

These two sets are equal up to the sign. Here,

$$\Phi(0) = \frac{1+\sqrt{5}}{2} = 1.6180339887, \quad \Phi(\pi) = \frac{-1+\sqrt{5}}{2} = 0.6180339887.$$

```
F0=golden_ratio1(0); fprintf('Phi(%8.6f) equals to %12.10f \n',0,F0);
>> Phi(0.000000) equals to 1.6180339887
F1=golden_ratio1(pi); fprintf('Phi(%8.6f) equals to %12.10f \n',pi,F1);
>> Phi(3.141593) equals to 0.6180339887
```

Note that $\Phi(0)\Phi(\pi) = 1$ and $\Phi(0) - \Phi(\pi) = 1$. Thus, for number $a>0$, the set of similarities equals to

$$S(a) = \{a\} = \{|a|\Phi(0), -|a|\Phi(\pi)\} = |a|\{\Phi(0), -\Phi(\pi)\} = |a|\{\Phi(0), 1 - \Phi(0)\}. \tag{22}$$

For the unit vector $a = e = 1$, the set of similarities is $S(e) = \{e\} = \{\Phi(0), 1 - \Phi(0)\}$. The golden pairs are $\{1, \Phi(0)\}$ and $\{1, 1 - \Phi(0)\}$.

***Example 2*** *(2-D vectors):*
Consider two vectors $\boldsymbol{a} = [a_1, a_2]$ and $\boldsymbol{b} = [b_1, b_2]$ in the 2D real space $R^2$. The inner product is defined as $(\boldsymbol{a}, \boldsymbol{b}) = a_1 b_1 + a_2 b_2$ and the norm of the vector $\boldsymbol{a}$ equals to $\|\boldsymbol{a}\| = \sqrt{a_1^2 + a_2^2}$. All unit vectors $\boldsymbol{e}$ have tips on the unit circle, that is, they are described by the set

$$\boldsymbol{e} \in \{\boldsymbol{e}_\varphi = (\cos \varphi, \sin \varphi); \; \varphi \in [0, 2\pi)\}.$$

We consider the polar form of the vector $\boldsymbol{a} = \|\boldsymbol{a}\|(\cos \alpha, \sin \alpha), \alpha \in [0, 2\pi)$, and the unit vector $\boldsymbol{e} = \boldsymbol{e}(\varphi)$ at angle $\varphi$ to the horizontal (see Fig. 11).

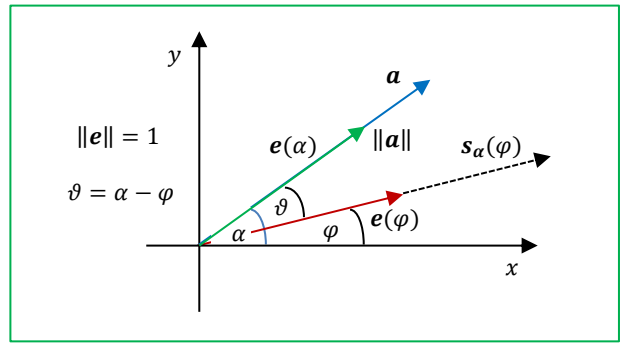

**Fig. 11** 2-D vectors.

Consider the unit vector $\boldsymbol{e}(\varphi)$ that composes the angle $\vartheta$ with the vector $\boldsymbol{a}$. The inner product of $\boldsymbol{e}(\varphi)$ with the unit vector $\boldsymbol{e}(\alpha) = \boldsymbol{a}/\|\boldsymbol{a}\|$ along the vector $\boldsymbol{a}$ equals to

$$\big(\boldsymbol{e}(\alpha), \boldsymbol{e}(\varphi)\big) = \cos \alpha \cos \varphi + \sin \alpha \sin \varphi = \cos(\alpha - \varphi) = \cos(\vartheta).$$

Along the angle $\varphi$, the vector that is in the golden ratio with the vector $\boldsymbol{a}$ equals to

$$\boldsymbol{s}_a(\varphi) = \|\boldsymbol{a}\|\Phi(\vartheta)\boldsymbol{e}(\varphi) = \|\boldsymbol{a}\|\Phi(\alpha - \varphi)\boldsymbol{e}(\varphi). \tag{23}$$

Therefore, the set of similarity of the vector $\boldsymbol{a}$ is defined as

$$S(\boldsymbol{a}) = \{\boldsymbol{s}_a(\varphi); \; \varphi \in [0, 2\pi)\} = \{\|\boldsymbol{a}\|\Phi(\alpha - \varphi)\boldsymbol{e}(\varphi); \; \varphi \in [0, 2\pi)\}. \tag{24}$$

Thus, for a given vector $\boldsymbol{a}$, a golden pair can be found along any direction. The golden ratio changes with angles. As an example, Fig. 12 shows the locus of all similar 2-D vectors $\boldsymbol{s}_a(\varphi)$ that compose the golden pairs $\{a, \boldsymbol{s}_a(\varphi)\}$, for the vectors $\boldsymbol{a} = [1, 2]$ and $[-1, 3]$ in parts (a) and (b), respectively. The similarity sets of these vectors are

$$S([1,2]) = \{\boldsymbol{s}_{[1,2]}(\varphi); \; \varphi \in [0, 2\pi)\} = \{\sqrt{5}\Phi(\tan^{-1} 2 - \varphi)\boldsymbol{e}(\varphi); \; \varphi \in [0, 2\pi)\}$$

and

$$S([-1,3]) = \{\boldsymbol{s}_{[-1,3]}(\varphi); \; \varphi \in [0, 2\pi)\} = \{\sqrt{10}\Phi(\pi - \tan^{-1} 3 - \varphi)\boldsymbol{e}(\varphi); \; \varphi \in [0, 2\pi)\}.$$

In these figures, the vectors are shown only for 128 uniformly distributed angles $\varphi$ from the interval $[0, 2\pi)$. The tips of the vectors are not shown. The figures recall the same petal rotated by different angles and magnifications.

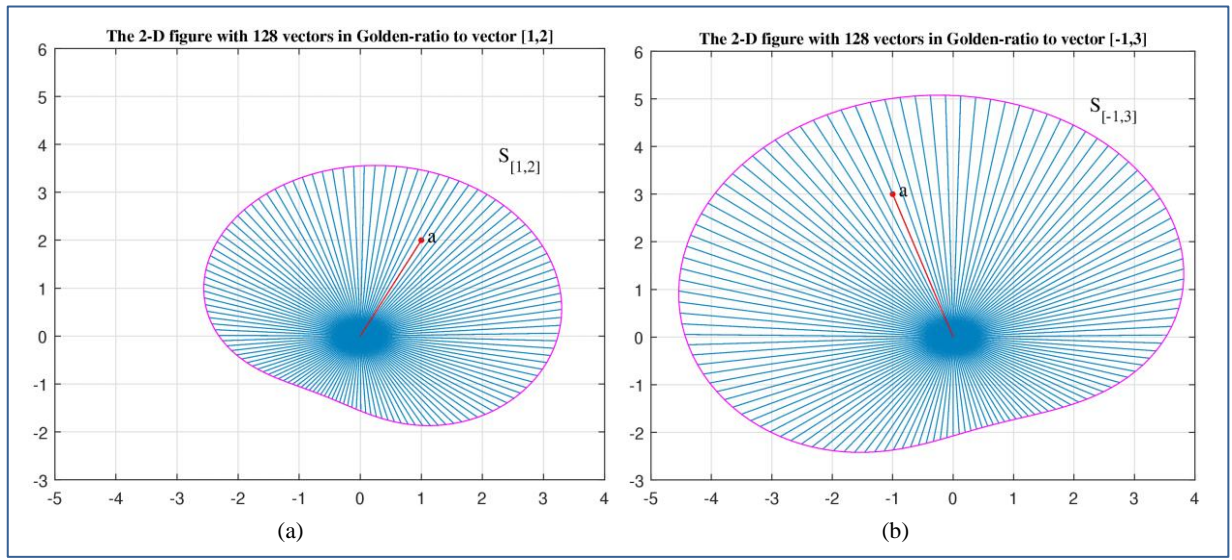

**Fig. 12** The locus of 128 golden pairs with the vectors (a) $a = [1,2]$ and (b) $a = [-1,3]$.

In part (c), the similar figure is shown for the vector $a = e = [1,0]$. The similarity set of this unit vector is

$$S([1,0]) = \{s_{[1,0]}(\varphi); \varphi \in [0,2\pi)\} = \{\Phi(-\varphi)e(\varphi); \varphi \in [0,2\pi)\} = \{\Phi(\varphi)e(\varphi); \varphi \in [0,2\pi)\}.$$

It should be noted that the figures of the similar sets in part (a) and (b) are the rotated figures of part (c) with magnification by the norms $\sqrt{5}$ and $\sqrt{10}$ of the vectors $a = [1,2]$ and $[-1,3]$, respectively. The figure for the set $S([1,2])$ can be obtained from the figure of the set $S([1,0])$, by the angle $\tan^{-1} 2$ and magnified by the number $\sqrt{5}$. Similarly, The figure of the similar set $S([1,2])$ in part (b) can be obtained from the figure of $S([1,0])$ by the angle $(\pi - \tan^{-1} 3)$ and magnified by the number $\sqrt{10}$.

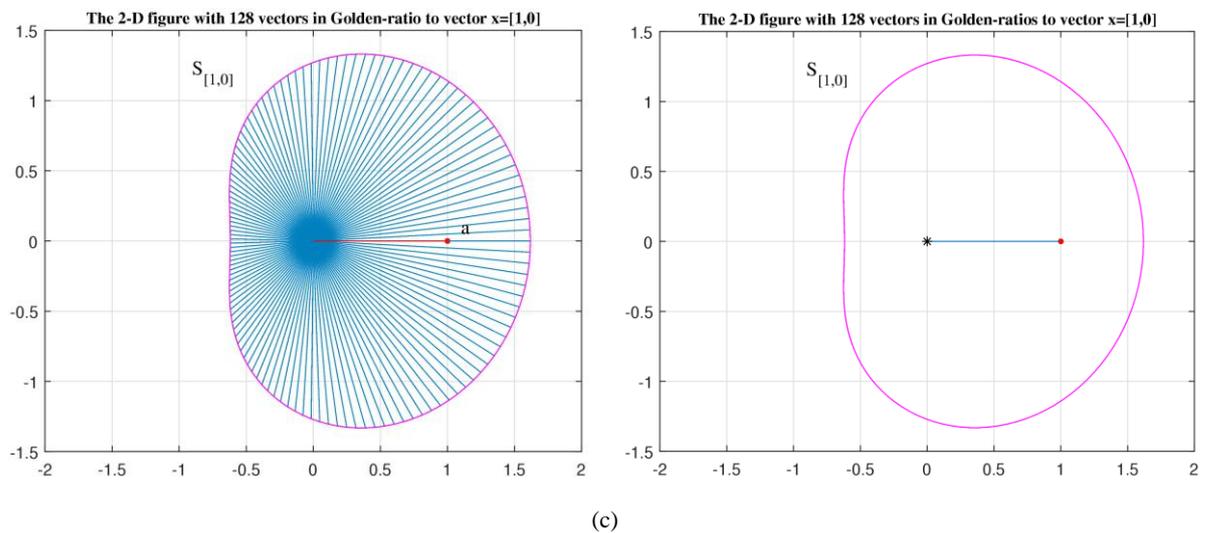

(c)

**Fig. 12** (continuation) (c) The locus of 128 Golden pairs with the unit vector $e = [1,0]$.

## 5. Field of Similarities

**A.** *Philosophical digression:* What is similarity in our case? Each vector affects its environment, stimulating its influence through the imposition of its likeness. If you think about it, then we are all a certain vector of possibilities that we impose on the environment by projecting ourselves into it, and this projection is symbolically represented by a certain projection angle. ("All our immediate environment is our projection, our likeness. This is a kind of vector shadow that is cast on the environment, and the Sun's ray symbolizes the angle of objectification of this shadow.") All of this needs to be well thought out

**B.** Until now, we knew that if two vectors do not interact and the inner product is zero (that is, the angle is 90 degrees), then mutual influence is excluded. But what is interesting is that the similarity coefficient at this angle is not equal to zero! That is, the influence is still there. According to the printout, it appears that

$$\Phi\left(\frac{\pi}{2}\right) = \sqrt{\Phi(0)} = \sqrt{\frac{1+\sqrt{5}}{2}} = 1.2720196495 \quad \text{(verified)}. \tag{25}$$

```
F2=golden_ratio1(pi/2); fprintf('Phi(%8.6f) equals to %12.10f \n',pi/2,F2);
>> Phi(1.570796) equals to 1.2720196495
F=sqrt(F0); fprintf('sqrt(%8.6f) equals to %12.10f \n',F0,F);
>> sqrt(1.618034) equals to 1.2720196495
```

**C.** Two roots of Eq. 7 are real, $\Phi_1(\alpha)$ and $\Phi_2(\alpha)$. We can say that being positive or negative number are two states, like heads and tails in the probability theory. After all, it was not for nothing that we got two states, one refers to positive roots, and the other to negative ones. Two other solutions are complex, $\Phi_3(\alpha)$ and $\Phi_4(\alpha)$. They show us our 2-D representation (the real solutions determine the 1D representation). Note that in English the words imaginary (for complex numbers) and the image have the same root.

**D. The sum of similarity vectors**

The following question arises. Is it possible to add similarity fields? Let us consider two different vectors $\boldsymbol{a}_1$ and $\boldsymbol{a}_2$ at angles of $\alpha_1 = \arg(\boldsymbol{a}_1)$ and $\alpha_2 = \arg(\boldsymbol{a}_2)$ to the positive real axis, respectively. The corresponding sets of similarities are

$$S(\boldsymbol{a}_1) = \{\boldsymbol{s}_{\boldsymbol{a}_1}(\varphi); \varphi \in [0,2\pi)\} = \{\|\boldsymbol{a}_1\|\Phi(\alpha_1 - \varphi)\boldsymbol{e}(\varphi); \varphi \in [0,2\pi)\}, \tag{26}$$

$$S(\boldsymbol{a}_2) = \{\boldsymbol{s}_{\boldsymbol{a}_2}(\varphi); \varphi \in [0,2\pi)\} = \{\|\boldsymbol{a}_2\|\Phi(\alpha_2 - \varphi)\boldsymbol{e}(\varphi); \varphi \in [0,2\pi)\}. \tag{27}$$

These two sets can be written as

$$S(\boldsymbol{a}_1) = \{\|\boldsymbol{a}_1\|\Phi(\varphi)\boldsymbol{e}(\alpha_1 - \varphi); \varphi \in [0,2\pi)\}$$
$$S(\boldsymbol{a}_2) = \{\|\boldsymbol{a}_2\|\Phi(\varphi)\boldsymbol{e}(\alpha_2 - \varphi); \varphi \in [0,2\pi)\}.$$

Then, their sum should be defined as the set of similarities

$$S(\boldsymbol{a}_1 + \boldsymbol{a}_2) = \{\|\boldsymbol{a}_1 + \boldsymbol{a}_2\|\Phi(\gamma - \varphi)\boldsymbol{e}(\varphi); \varphi \in [0,2\pi)\}, \tag{28}$$

where the angle $\gamma = \arg(\boldsymbol{a}_1 + \boldsymbol{a}_2)$. Let us verify if the following is true:

$$S(\boldsymbol{a}_1 + \boldsymbol{a}_2) = S(\boldsymbol{a}_1) + S(\boldsymbol{a}_2). \tag{29}$$

Here, the summation is angle-wise, that is, the summation of vectors that correspond to the same angle $\varphi$. Therefore, this equation can be written as

$$\|\boldsymbol{a}_1 + \boldsymbol{a}_2\|\Phi(\varphi)\boldsymbol{e}(\gamma - \varphi) = \|\boldsymbol{a}_1\|\Phi(\varphi)\boldsymbol{e}(\alpha_1 - \varphi) + \|\boldsymbol{a}_2\|\Phi(\varphi)\boldsymbol{e}(\alpha_2 - \varphi). \tag{30}$$

Here, all vectors have the same coefficient of similarity. Removing the similar term $\Phi(\varphi)$ from this equation, we obtain

$$\|\boldsymbol{a}_1 + \boldsymbol{a}_2\|\boldsymbol{e}(\gamma - \varphi) = \|\boldsymbol{a}_1\|\boldsymbol{e}(\alpha_1 - \varphi) + \|\boldsymbol{a}_2\|\boldsymbol{e}(\alpha_2 - \varphi). \tag{31}$$

This equation describes the well-known rule for summing vectors over projections. The following calculations are valid for the right part of this equation:

$$\|\boldsymbol{a}_1\|\boldsymbol{e}(\alpha_1 - \varphi) + \|\boldsymbol{a}_2\|\boldsymbol{e}(\alpha_2 - \varphi) = \|\boldsymbol{a}_1\|[\cos(\alpha_1 - \varphi), \sin(\alpha_1 - \varphi)] + \|\boldsymbol{a}_2\|[\cos(\alpha_2 - \varphi), \sin(\alpha_2 - \varphi)]$$
$$= \|\boldsymbol{a}_1\|[\cos(\alpha_1)\cos(\varphi) + \sin(\alpha_1)\sin(\varphi), \sin(\alpha_1)\cos(\varphi) - \cos(\alpha_1)\sin(\varphi)]$$
$$+ \|\boldsymbol{a}_2\|[\cos(\alpha_2)\cos(\varphi) + \sin(\alpha_2)\sin(\varphi), \sin(\alpha_2)\cos(\varphi) - \cos(\alpha_2)\sin(\varphi)]$$
$$= [(\boldsymbol{a}_1)_x \cos(\varphi) + (\boldsymbol{a}_1)_y \sin(\varphi), (\boldsymbol{a}_1)_y \cos(\varphi) - (\boldsymbol{a}_1)_x \sin(\varphi)]$$
$$+ [(\boldsymbol{a}_2)_x \cos(\varphi) + (\boldsymbol{a}_2)_y \sin(\varphi), (\boldsymbol{a}_2)_y \cos(\varphi) - (\boldsymbol{a}_2)_x \sin(\varphi)]$$
$$= [((\boldsymbol{a}_1)_x + (\boldsymbol{a}_2)_x)\cos(\varphi) + ((\boldsymbol{a}_1)_y + (\boldsymbol{a}_2)_y)\sin(\varphi), ((\boldsymbol{a}_1)_y + (\boldsymbol{a}_2)_y)\cos(\varphi)$$
$$- ((\boldsymbol{a}_1)_x + (\boldsymbol{a}_2)_x)\sin(\varphi)]$$
$$= [(\boldsymbol{a}_1 + \boldsymbol{a}_2)_x \cos(\varphi) + (\boldsymbol{a}_1 + \boldsymbol{a}_2)_y \sin(\varphi), (\boldsymbol{a}_1 + \boldsymbol{a}_2)_y \cos(\varphi) - (\boldsymbol{a}_1 + \boldsymbol{a}_2)_x \sin(\varphi)]$$
$$= \|\boldsymbol{a}_1 + \boldsymbol{a}_2\|[\cos(\gamma - \varphi), \sin(\gamma - \varphi)] = \|\boldsymbol{a}_1 + \boldsymbol{a}_2\|\boldsymbol{e}(\gamma - \varphi).$$

Here, $(\boldsymbol{a}_1)_x = \|\boldsymbol{a}_1\|\cos(\alpha_1)$, $(\boldsymbol{a}_1)_y = \|\boldsymbol{a}_1\|\sin(\alpha_1)$, $(\boldsymbol{a}_2)_x = \|\boldsymbol{a}_2\|\cos(\alpha_2)$, $(\boldsymbol{a}_2)_y = \|\boldsymbol{a}_2\|\sin(\alpha_2)$, and $(\boldsymbol{a}_1 + \boldsymbol{a}_2)_x = \|\boldsymbol{a}_1 + \boldsymbol{a}_2\|\cos(\gamma)$, $(\boldsymbol{a}_1 + \boldsymbol{a}_2)_y = \|\boldsymbol{a}_1 + \boldsymbol{a}_2\|\sin(\gamma)$. Thus, everything is correct in Eqs. 30 and 31.

Figure 13 illustrates this property in part (b), where the parallelogram is the result of the rotation of the original parallelogram in part (a), which is composed for the sum of vectors $\boldsymbol{a}_1$ and $\boldsymbol{a}_2$.

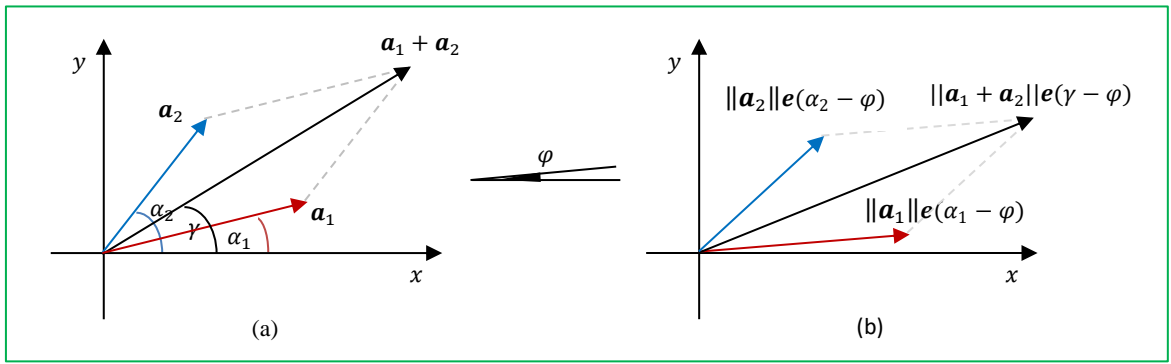

**Fig. 13** The parallelogram composed of the sum of two vectors (a) before and (b) after the rotation by the angle $\varphi$.

As examples, Figure 14 illustrates the summation of the similarity sets for the vectors $\boldsymbol{a}_1 = [1,2]$ and $\boldsymbol{a}_2 = [-1,3]$ in part (a) and for vectors $\boldsymbol{a}_1 = [2,-3]$ and $\boldsymbol{a}_2 = [1,5]$.

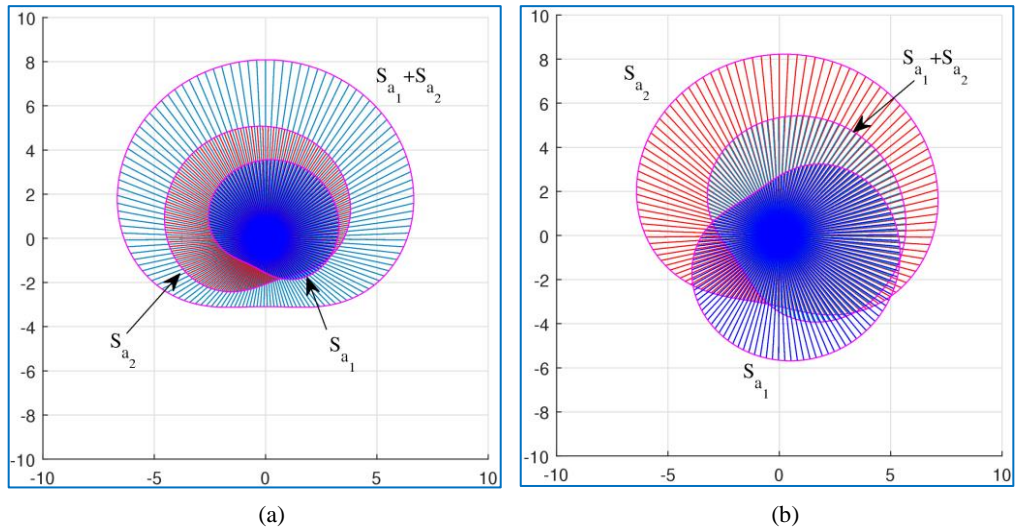

**Fig. 14** Sum of similarity sets for the 2-D vectors (a) $\boldsymbol{a}_1 = [1,2]$, $\boldsymbol{a}_2 = [-1,3]$ and (b) $\boldsymbol{a}_1 = [2,-3]$, $\boldsymbol{a}_2 = [1,5]$.

Figure 15 shows the same sum $S([1,1])$ of similarity sets, for three pairs of vectors $\boldsymbol{a}_1$ and $\boldsymbol{a}_2$. Namely, for the vectors $[0,1]$ and $[1,0]$ in part (a), vectors $[1,2]$ and $[0,-1]$ in part (b), and vectors $[2,3]$ and $[-1,-2]$ in part (c).

*Art, January 12, 2025*

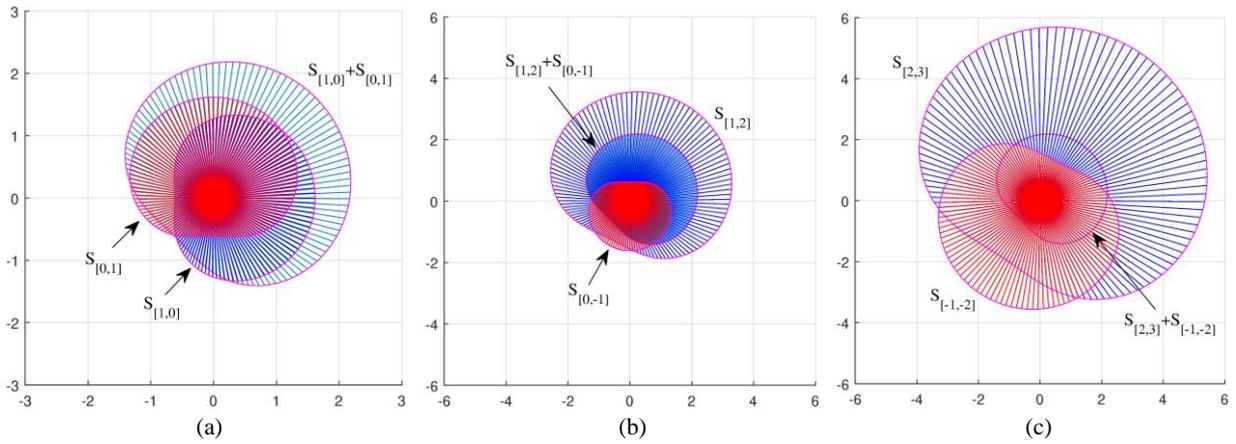

**Fig. 15** Sum of similarity sets for the 2D vectors (a) $a_1 = [0,1]$, $a_2 = [1,0]$, (b) $a_1 = [1,2]$, $a_2 = [0,-1]$, and (c) $a_1 = [2,3]$, $a_2 = [-1,-2]$.

It should be noted that we do not sum the similarity sets over equally directed rays. If we do that, that is, if consider the sum of the similar vectors

$$A = \|a_1\|\Phi(\alpha_1 - \varphi)e(\varphi) + \|a_2\|\Phi(\alpha_2 - \varphi)e(\varphi)$$

for each angle $\varphi \in [0, 2\pi]$, then, we need to find the vector $a_0$ and angle $\gamma_0$, such that $\|a_0\|\Phi(\gamma_0 - \varphi)e(\varphi) = A$. The solution of the equation

$$\|a_0\|\Phi(\gamma_0 - \varphi) = \|a_1\|\Phi(\alpha_1 - \varphi) + \|a_2\|\Phi(\alpha_2 - \varphi), \qquad \varphi \in [0, 2\pi),$$

is unknown for us.

***Example 3*** *(3-D vectors):* We consider the traditional representation of the 3D unit vectors, namely the set

$$\{e\} = \{[\sin\varphi \cos\psi, \sin\varphi \sin\psi, \cos\varphi]; \varphi \in [0,\pi], \psi \in [0,2\pi)\}. \tag{32}$$

The geometry of such a unit vector $e = e_{\varphi,\psi} = [\sin\varphi \cos\psi, \sin\varphi \sin\psi, \cos\varphi]$ is shown in Fig. 16.

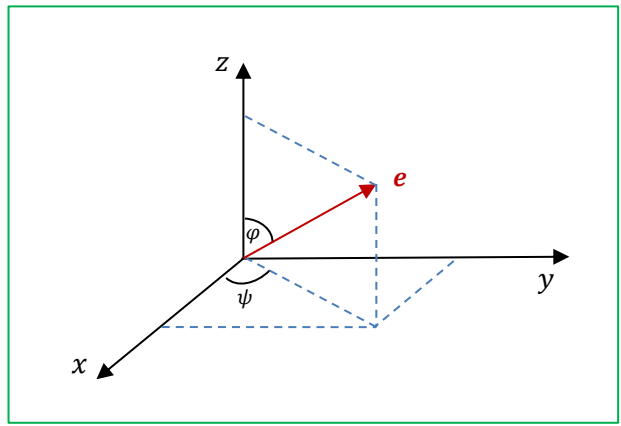

**Fig. 16** The unit vector $e$ in the 3D spherical coordinate system.

The inner product of such a unit vector $e$ with a vector $a = \|a\|[\sin\varphi_1 \cos\psi_1, \sin\varphi_1 \sin\psi_1, \cos\varphi_1]$, where $\varphi_1 \in [0,\pi], \psi_1 \in [0,2\pi)$, is calculated by

$$(e, a) = \|a\|\cos\vartheta = \|a\|[\sin\varphi \cos\psi \sin\varphi_1 \cos\psi_1 + \sin\varphi \sin\psi \sin\varphi_1 \sin\psi_1 + \cos\varphi \cos\varphi_1].$$

Thus,
$$\cos\vartheta = \sin\varphi \cos\psi \sin\varphi_1 \cos\psi_1 + \sin\varphi \sin\psi \sin\varphi_1 \sin\psi_1 + \cos\varphi \cos\varphi_1. \tag{33}$$

The cosine of the angle between these two vectors is the function of 4 arguments, that is, the angle $\vartheta$ between the vectors $\boldsymbol{e}$ and $\boldsymbol{a}$ is the function $\vartheta = \vartheta(\varphi, \psi, \varphi_1, \psi_1)$. The similarity set, that is, the set of all vectors in golden ratios with the vector $\boldsymbol{a}$ is

$$S(\boldsymbol{a}) = \{\|\boldsymbol{a}\|\Phi(\vartheta)[\sin\varphi\cos\psi, \sin\varphi\sin\psi, \cos\varphi]; \varphi \in [0,\pi], \psi \in [0,2\pi]\}. \tag{34}$$

Let $\boldsymbol{a}$ be the unit vector be $\boldsymbol{e}_3 = [0,0,1]$. Figure 17 shows the locus of vectors in Golden ratio with this vector in part (a). In this case, $\varphi_1 = 0$ and $\cos\vartheta = \cos\vartheta(\varphi, \psi, 0, \psi_1) = \cos\varphi$. Therefore, $\vartheta = \varphi$ and the set

$$S(\boldsymbol{e}_3) = \{\Phi(\varphi)[\sin\varphi\cos\psi, \sin\varphi\sin\psi, \cos\varphi]; \varphi \in [0,\pi], \psi \in [0,2\pi]\}. \tag{35}$$

We also consider the unit vector $\boldsymbol{a} = \boldsymbol{e}_1 = [1,0,0]$. Then, $\varphi_1 = \pi/2$ and $\psi_1 = 0$ and $\cos\vartheta = \cos\vartheta(\varphi, \psi, \pi/2, 0) = \sin\varphi\cos\psi$. Therefore, $\vartheta = \arccos(\sin\varphi\cos\psi)$ and the similar set is

$$S(\boldsymbol{e}_1) = \{\Phi(\vartheta)[\sin\varphi\cos\psi, \sin\varphi\sin\psi, \cos\varphi]; \varphi \in [0,\pi], \psi \in [0,2\pi]\}. \tag{36}$$

The locus of vectors of this similar set is shown in part (b). In these two figures, as in the 2D case above, only the ends of the vectors as dots are shown. The angles $\varphi$ and $\psi$ were taken with the step $2\pi/512$ in the intervals $[0,\pi]$ and $[0,2\pi]$, respectively.

Figures 18 and 19 show these similar sets by different angles. Namely, by using the azimuth (AZ) of zero degree and vertical elevation (EL) of 90 and 180 degrees, respectively. For this, the MATLAB functions 'view(2)' and 'view(0,180)' were used.

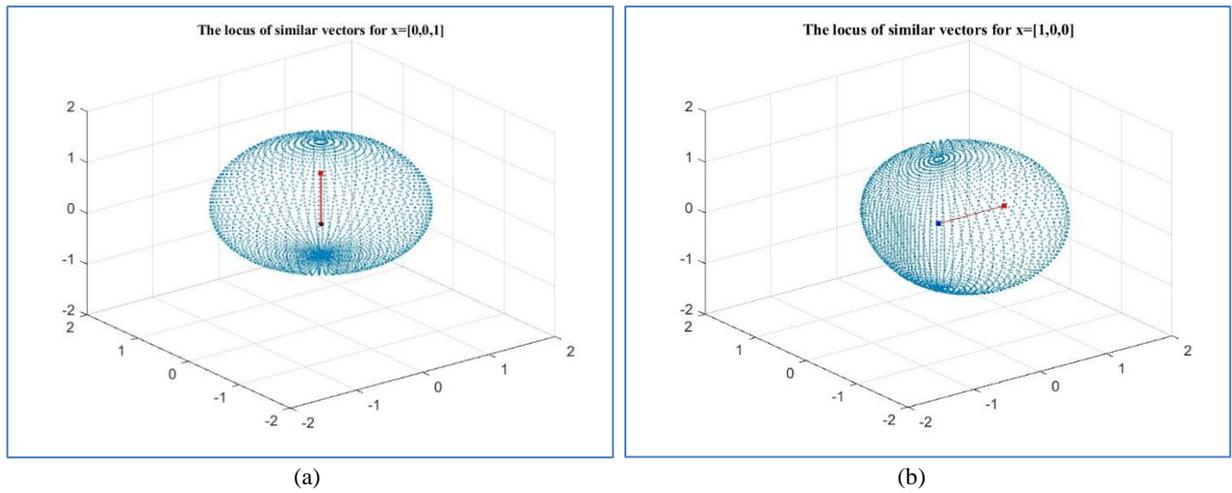

(a)  (b)

**Fig. 17** The geometry of the similar sets of golden vectors with the 3D vectors (a) [0,0,1] and (b) [1,0,0].

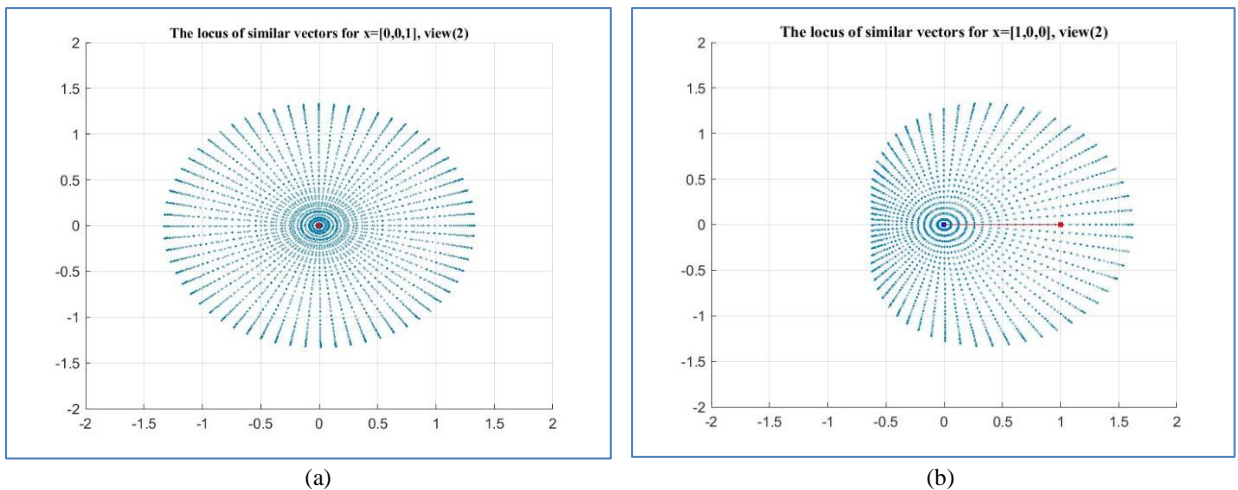

(a)  (b)

**Fig. 18** The view with AZ=0 and EL=90 degrees of geometry of the similar sets of golden vectors with the 3D vectors (a) [0,0,1] and (b) [1,0,0].

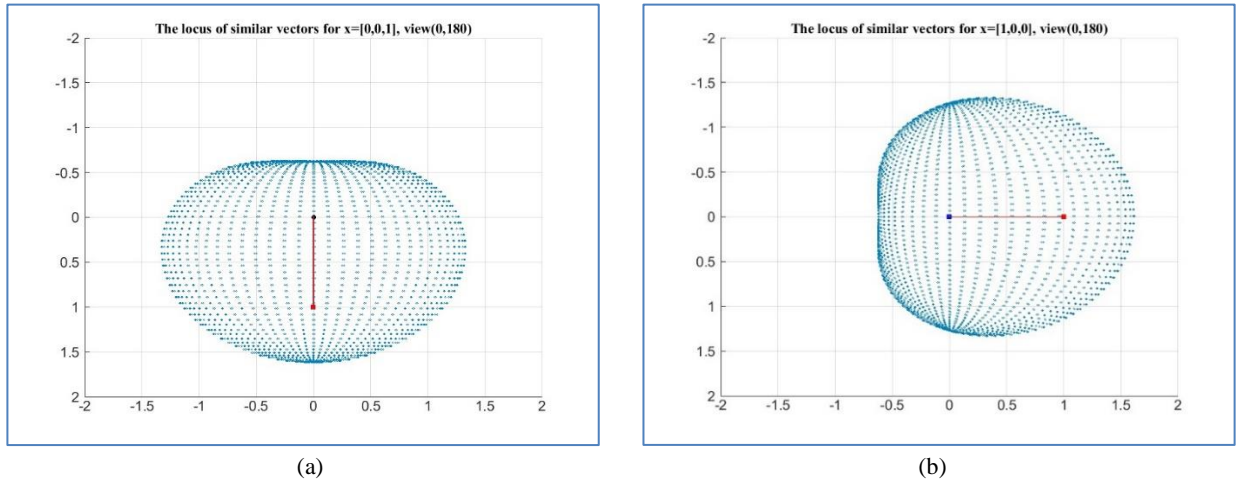

**Fig. 19** The view with AZ=0 and EL=180 degrees: The geometry of the similar sets of golden vectors with the 3-D vectors (a) [0,0,1] and (b) [1,0,0].

Figure 20 shows the 3D surface that is made up of the vertices of the vectors of a subset of $S(e_3)$ in part (a). Thus, this is the surface framing the vectors which are similar to the unit vector $e_3 = [0,0,1]$. In part (b), the similar surface is shown for the unit vector $e_1 = [1,0,0]$. For both surfaces, the angles $\varphi$ and $\psi$ were taken with the step $2\pi/360$ in the intervals $[0, \pi]$ and $[0, 2\pi]$, respectively.

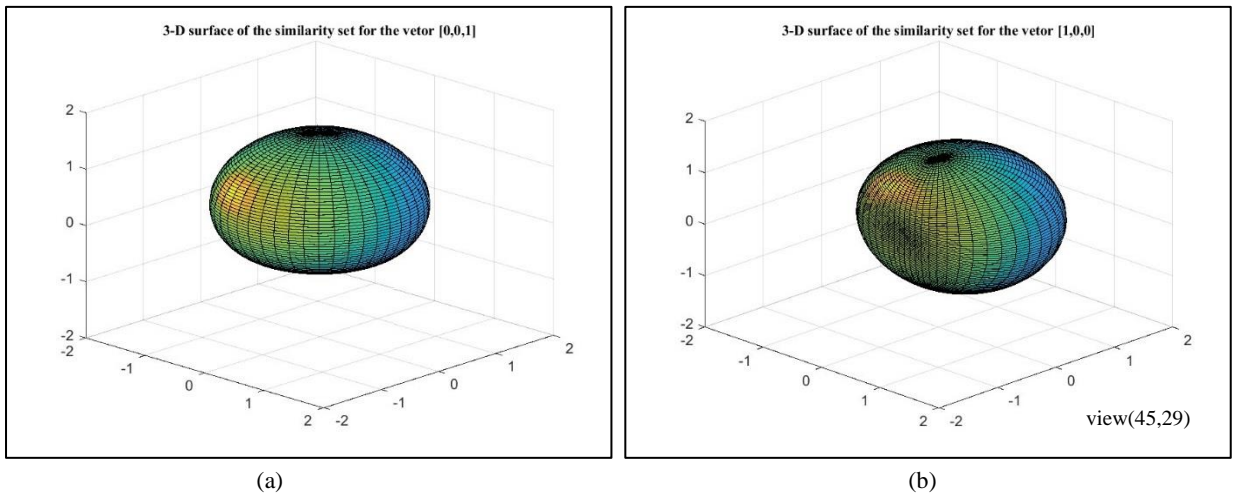

**Fig. 20** The 3-D surface of the similar sets of golden vectors with the 3-D vectors (a) [0,0,1] and (b) [1,0,0].

## 6. Similarity of Triangles

In this section, we consider triangles as elements of the 6D vector space and introduce the concept of the inner product and norm of triangles. The tringles in golden ratio are described and the similarity sets are presented with examples.

In order to show the similarity of three points $(a, b, c)$ in the form of a triangle on the plane (see Fig. 21(a)), we need three 2D vectors, which we denote by $v_a, v_b$, and $v_c$. The vectors $v_a = (x_a, y_a)$, $v_b = (x_b, y_b)$, and $v_c = (x_c, y_c)$. These coordinate vectors compose the 6-D vector $V = (v_a, v_b, v_c)$.

Consider two 6D vectors that correspond to 2 triangles, $V_1 = (v_{a_1}, v_{b_1}, v_{c_1})$ and $V_2 = (v_{a_2}, v_{b_2}, v_{c_2})$. The inner product of these vectors can be defined as

$$(V_1, V_2) = (v_{a_1} - v_{b_1}, v_{a_2} - v_{b_2}) + (v_{b_1} - v_{c_1}, v_{b_2} - v_{c_2}) + (v_{c_1} - v_{a_1}, v_{c_2} - v_{a_2}). \tag{37}$$

The norm of the vector is defined as

$$\|V\|^2 = (V, V) = \|v_a - v_b\|^2 + \|v_b - v_c\|^2 + \|v_c - v_a\|^2. \tag{38}$$

The norm $\|V\| = 0$, when a triangle degenerates into point, that is, when $v_a = v_b = v_c$, and this case is not considered.

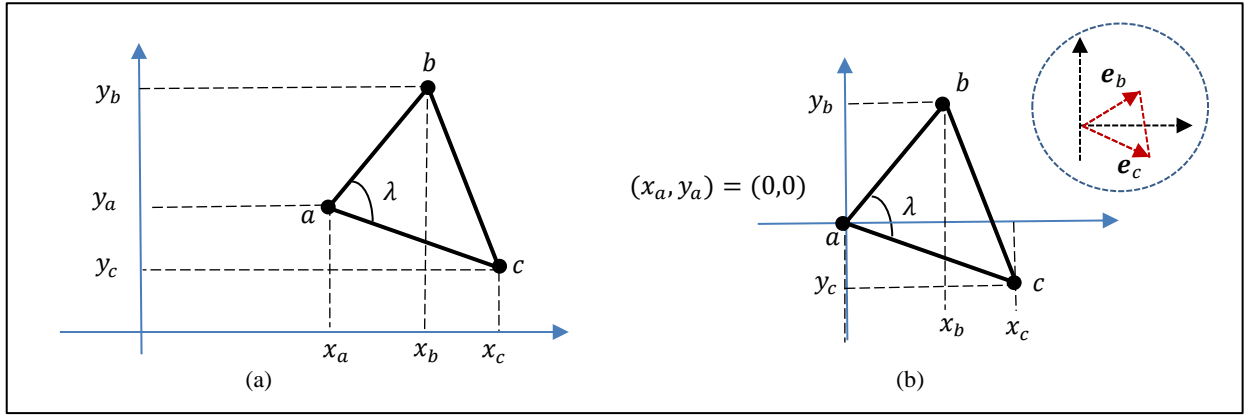

**Fig. 21** (a) Triangle for the vector $V$ and (b) the triangle with the vertex $a$ at the origin.

A unit vector, or a triangle, $E = (e_a, e_b, e_c)$ with norm 1 is defined as

$$\|E\|^2 = \|e_a - e_b\|^2 + \|e_b - e_c\|^2 + \|e_c - e_a\|^2 = 1. \tag{39}$$

We can zero the first 2-D vector $e_a$ and consider the unit vector $E = (0, e_b, e_c)$, for which

$$\|E\|^2 = \|e_b\|^2 + \|e_b - e_c\|^2 + \|e_c\|^2 = 1, \tag{40}$$

with condition that $e_b \neq e_c$, $\|e_b\| \neq 0$, and $\|e_c\| \neq 0$. Such 6D vector corresponds to a triangle with the point $a$ at the center of the system of coordinates. Such an example is shown in Fig. 21 in part (b).

It is follows from Eq. 40, that

$$2\|e_b\|^2 - 2\|e_b\|\|e_c\|\cos(\lambda) + 2\|e_c\|^2 = 1.$$

Here, $\lambda$ is the angle between 2D vectors $e_b$ and $e_c$. This equation can be written as

$$\left(\|e_b\| - \frac{1}{2}\|e_c\|\cos(\lambda)\right)^2 + \frac{1}{4}\|e_c\|^2[4 - \cos^2(\lambda)] = \frac{1}{2}. \tag{41}$$

Solutions of this equation can be written as ($\lambda \neq 0$)

$$\begin{cases} \|e_b\| - \frac{1}{2}\|e_c\|\cos(\lambda) = \frac{1}{\sqrt{2}}\cos(\phi), \\ \frac{1}{2}\|e_c\|\sqrt{4 - \cos^2(\lambda)} = \frac{1}{\sqrt{2}}\sin(\phi). \end{cases} \tag{42}$$

Thus, we have a parameterized set of solutions; two parameters are the angles $\lambda$ and $\phi$. The solutions can be written as

$$\begin{cases} \|e_c\| = \sqrt{2}/\sqrt{4 - \cos^2(\lambda)}\sin(\phi), \\ \|e_b\| = \frac{1}{2}\|e_c\|\cos(\lambda) + \frac{1}{\sqrt{2}}\cos(\phi). \end{cases} \tag{43}$$

Here, $0 < \phi \leq \pi - \lambda$ and $0 < \lambda < \pi$. This system of solutions connects the lengths and the angle $\lambda$ between the sides of the triangle, $e_b$ and $e_c$. Note that, to generalize this solution, we can add a zero element $V_c = (c, c, c)$ with norm 0. Indeed, for a 2D vector $c = (c_1, c_2) \neq (0,0)$, $\|E + V_c\| = \|E\| = 1$.

The vector $E = (e_a, e_b, e_c)$ can be analytically written as

$$\begin{cases} e_a = 0, \\ e_c = \|e_c\|[\cos(\phi), \sin(\phi)], \\ e_b = \|e_b\|[\cos(\lambda + \phi), \sin(\lambda + \phi)]. \end{cases} \tag{44}$$

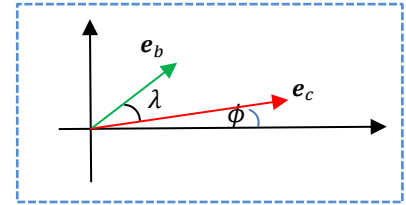

The angles $\phi$ in Eqs. 43 and 44 are considered the same.

Thus, the vector is parameterized by two angles, that is, $\boldsymbol{E} = \boldsymbol{E}(\phi, \lambda)$. In this system, the vector $\boldsymbol{e}_c$ is rotated counter clock-wise by the angle $\phi$ and the vector $\boldsymbol{e}_c$ by the angle $(\lambda + \phi)$ (from the horizontal). We denote the set of such unit 6D vectors (triangles) $\boldsymbol{E}$ by $\mathcal{E}$. As examples, Fig. 22 shows five unit triangles with the angle $\lambda = 40^o$. The first triangle with angle $\phi = 10^o$ is shown in red with vertices marked. Other four unit triangles are shown with the angles $\phi = 40^o, 70^o, 100^o$, and $130^o$.

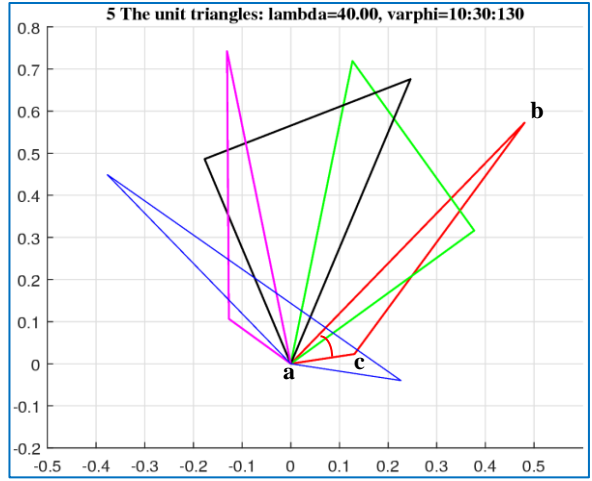

**Fig. 22** Five unit triangles with the angle $\lambda = 40^o$.

Two other examples with unit triangles are shown in Fig. 23. Four unit triangles with angle $\lambda = 60^o$ are shown in part (a), when angles $\phi = 10^o, 40^o, 70^o$, and $100^o$. In part (b), three unit triangles with angle $\lambda = 110^o$ are shown ), when angles $\phi = 10^o, 40^o$, and $70^o$.

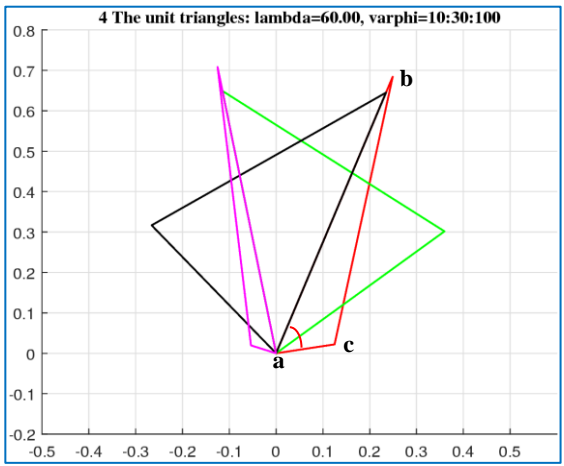
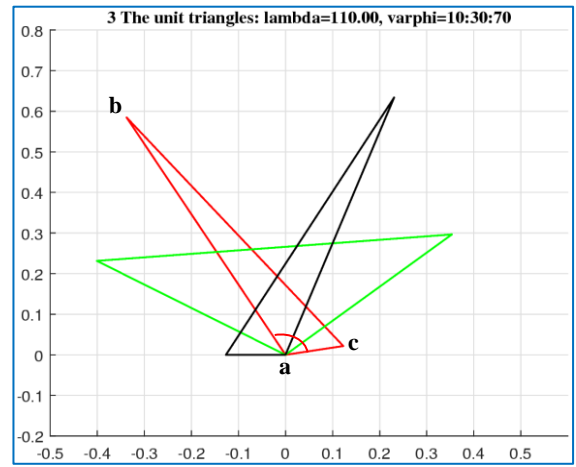

**Fig. 23** Unit triangles with the angle (a) $\lambda = 60^o$ and (b) $\lambda = 110^o$.

For a 6D vector $\boldsymbol{V} = (\boldsymbol{v}_a, \boldsymbol{v}_b, \boldsymbol{v}_c)$ presenting a triangle, the similar triangle, or the triangle in the general golden ratio with $\boldsymbol{V}$, is defined as

$$s_V(\phi, \lambda) = \|\boldsymbol{V}\|\Phi(\vartheta)\boldsymbol{E}(\phi, \lambda), \quad 0 < \phi \le \pi - \lambda, \lambda \in (0, \pi). \qquad (45)$$

Here, $\vartheta = angle(\boldsymbol{V}, \boldsymbol{E})$ denotes the angle between vectors $\boldsymbol{V}$ and $\boldsymbol{E} = (\boldsymbol{e}_a, \boldsymbol{e}_b, \boldsymbol{e}_c)$, which is calculated by

$$\cos \vartheta = \frac{(\boldsymbol{V}, \boldsymbol{E})}{\|\boldsymbol{V}\|}.$$

The norm $\|\boldsymbol{V}\|$ is calculated from Eq. 39 and the inner product is calculated by

$$(V, E) = (v_a - v_b, e_a - e_b) + (v_b - v_c, e_b - e_c) + (v_c - v_a, e_c - e_a). \tag{46}$$

As an example, Figure 24 show the triangle described by the vector $V = (v_a, v_b, v_c) = ([0,0], [3,2], [5,0])$ in part (a). The angle between the vectors $v_b$ and $v_c$ equals to $\theta = 33.69^o$ and $\Phi(\vartheta) = 1.5702$. The similarity triangles $s_V(20^o, \lambda)$, $s_V(40^o, \lambda)$, and $s_V(70^o, \lambda)$ for angle $\lambda = \theta$, are shown in parts (b), (c), and (d), respectively.

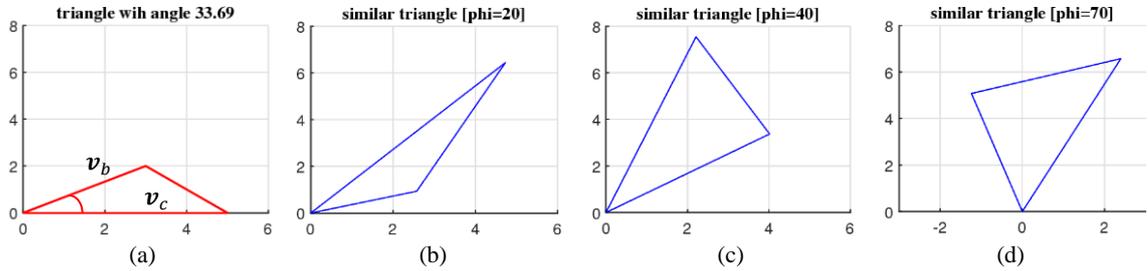

**Fig. 24** (a) The original triangle with angle $\theta = 33.69^o$ and similarity triangles with the angle $\phi$ of (b) $20^o$, (c) $40^o$, and (d) $70^o$.

For the vector $V = (v_a, v_b, v_c)$, the similarity set of triangles is defined as

$$S(V) = \{s_V(\phi, \lambda); \ 0 < \phi \le \pi - \lambda, 0 < \lambda < \pi\}. \tag{47}$$

Also, we can write this set as $S(V) = \{\|V\|\Phi(\vartheta)E; \ E \in \mathcal{E}, \vartheta = angle(V, E)\}$.

As examples, Figure 25 in part (a) shows the subset of the similarity set $S(V)$, for the vector $V = ([0,0], [0,2], [3,2])$ which represents a right triangle with angle $\theta = 56.31^o$ between sides ab and ac (shown in red). The angle $\theta = \text{angle}(v_b, v_c)$ is the angle between vectors $v_b$ and $v_c$. The unit vectors $E$ are calculated by Eq. 43 and 44 for 23 angles $\phi = 10^o: 5^o: 120^o$. The second parameter $\lambda = \theta$, that is, the angle between the vectors $e_b$ and $e_c$ in the similar triangles is the same as in the triangle for $V$. In part (b), the subset of another similarity set $S(V)$ is shown, for the vector $V = ([0,0], [0,2], [3,3])$. This vector represents a triangle with angle of $\theta = 45^o$ between vectors $v_b$ and $v_c$ (shown in red). 26 angles of $\phi$ are used, namely $\phi = 10^o: 5^o: 135^o$, and the angle $\lambda = \theta$.

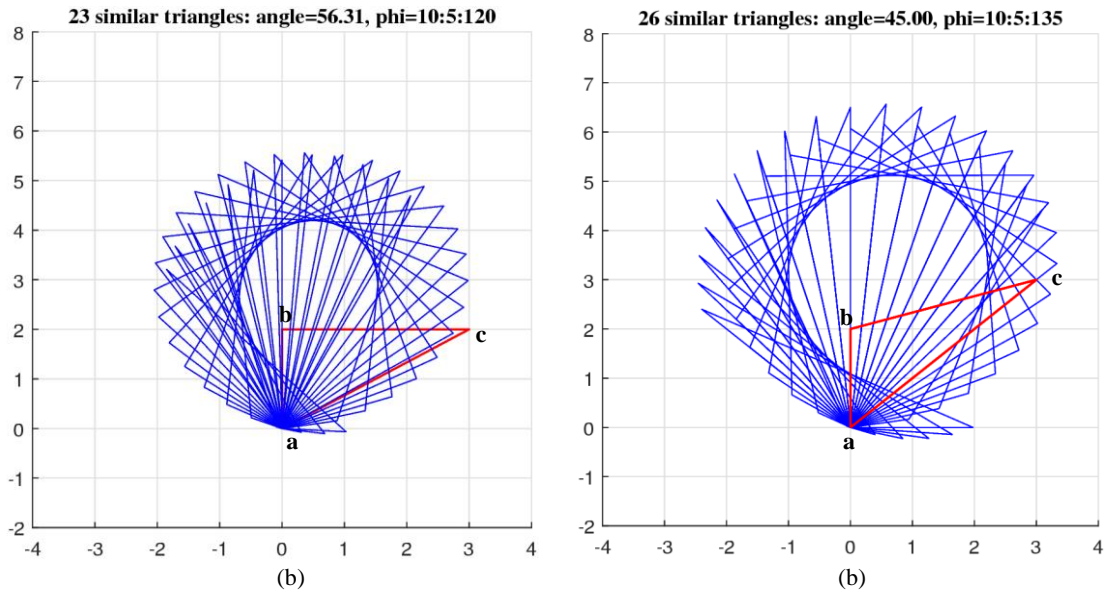

**Fig. 25** The similar triangles for (a) the triangle with angle $\theta = 56.31^o$ and (b) the triangle with angle $\theta = 45^o$.

In Figure 26, the similar subsets are shown for the equilateral triangle with sides of length 5. The corresponding vector is $V = ([1,2], [3.5, 4\sin(\pi/3)], [7,2])$. The first point $a$ of the triangle is not at the origin. In part (a), the subset of similar triangles is shown for 29 angles $\phi = 10^o: 5^o: 150^o$ and angle $\lambda = 60^o$. This is the case when $\lambda = \theta = 60^o$. In part (b), the subset is shown for the same 29 angles of $\phi$ and the angle $\lambda = 30^o$.

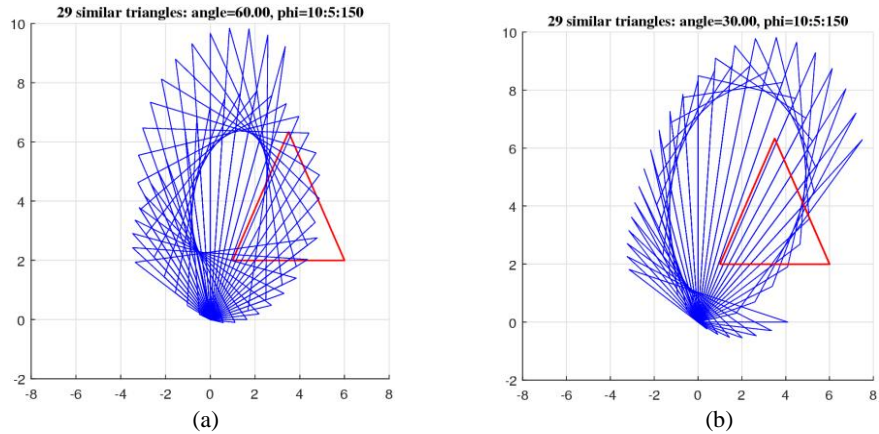

**Fig. 26** The similar triangles for the equilateral triangle when the angle (a) $\lambda = 60^o$ and (b) $\lambda = 30^o$.

What can note from the figures above that among similar triangles there is no equal to the triangle with the vector $\boldsymbol{V}$. It also not difficult to see from Eq. 45, that $s_V(\phi, \lambda) = \boldsymbol{V}$ only if

$$\Phi(\vartheta)\boldsymbol{E}(\phi, \lambda) = \frac{\boldsymbol{V}}{\|\boldsymbol{V}\|},$$

that is, when $\Phi(\vartheta) = 1$, or $\vartheta = 120^o$ and $\lambda = \vartheta$. It is possible that other definitions of the inner product of triangles could lead to similarities that include the original triangle $\boldsymbol{V}$.

As was mentioned above, to generalize this solution, we can add to unit vectors $\boldsymbol{E} = \boldsymbol{E}(\phi, \lambda)$ a zero element $\boldsymbol{V_c} = (\boldsymbol{c}, \boldsymbol{c}, \boldsymbol{c})$ with a 2-D vector $\boldsymbol{c} = (c_1, c_2) \neq (0,0)$. Therefore in general, we can consider the similarity set of triangles of the vector-triangle $\boldsymbol{V} = (\boldsymbol{v_a}, \boldsymbol{v_b}, \boldsymbol{v_c})$ as a set parameterized by angles $\phi, \lambda$, and vector $\boldsymbol{c}$,

$$S(\boldsymbol{V}) = \{s_V(\phi, \lambda, \boldsymbol{c}); \ 0 < \phi \leq \pi - \lambda, 0 < \lambda < \pi, \boldsymbol{c} \in R^2\}. \tag{48}$$

To facilitate understanding, we can separate similarities, by fixing two parameters out of three. Also, we can consider these similarities separately (as spatial similarity, similarity in one fixed angle $\lambda$, or similarity in rotation $\phi$). If we are interested in similarities with a fixed angle between two sides, then we should fix the value of $\lambda$ by giving it the value of one of the angles of the original triangle. Adding a constant vector $\boldsymbol{V_c}$ leads to a spatial shift (translation) of the set of similarities.

## 7. Afterword and Conclusion

In this work (Part I of our study), we have abstractly generalized the similarity law for multidimensional vectors. Initially, the law of similarity was derived for one-dimensional vectors. Although it operated with such values of the ratio of parts of the whole, it meant linear dimensions (a line is one-dimensionality). Now the question arises - where can one observe this generalization? If we refer to the forms of fruits, for instance the apple, then it is easy to see a strongly pronounced broken sphericity, which does correspond to the listed forms.

The main part of this research is the concept of similarity, namely the set of similar vectors to a given one. In all graphics above, the vectors were used, not area or surface (in the 3-D case). The vector represents the force. This concept is difficult to explain, even to understand. There is a mystery here that we can't solve. Intuitively our subconscious mind shows that a certain vector spreads itself in the surrounding space. What we consider local, as a certain center of outflow of a narrowly directing force, is an abstraction. Any force retains its likeness around itself and has spatial dimensions even in the opposite direction. Physically, it's hard to imagine. Psychologically, it can be imagined in the following way. If you accompany the realization of force, then you will have a gain ($\Phi(\alpha), \alpha = 0$), but your opposition to it is not completely destroyed ($\Phi(\alpha), \alpha = \pi$). Nature reserves the right to object with coefficient $\Phi(\pi) = (\sqrt{5} - 1)/2$, when there is a right to encourage with coefficient $\Phi(0) = (\sqrt{5} + 1)/2$. Also, $\Phi(0) \Phi(\pi) = 1$, which means that the impact of the external environment is intended only to separate those who are with it from those who are against it. Deviations of action along all other angles between standing and confrontation have a tendency to roll to one of these sides (the sign of the derivative of $\Phi(\alpha)$). In this perspective, the law of similarity is clear.